\documentclass{amsart}
\usepackage{graphicx}
\newtheorem{theorem}{Theorem}

\usepackage{amsmath,amscd}

\begin{document}

\title{A Topological Obstruction to the Removal of a Degenerate Complex Tangent and Some Related Homotopy and Homology Groups}

\author{Ali M. Elgindi}

\begin{abstract}
In this article, we derive a topological obstruction to the removal of a isolated degenerate complex tangent to an embedding of a 3-manifold into $\mathbb{C}^3$ (without affecting the structure of the remaining complex tangents). We demonstrate how the vanishing of this obstruction is both a necessary and sufficient condition for the (local) removal of the isolated complex tangent. The obstruction is a certain homotopy class of the space $\mathbb{Y}$ consisting of totally real 3-planes in the Grassmanian of real 3-planes in $\mathbb{C}^3$ (=$\mathbb{R}^6$). We further compute additional homotopy and homology groups for the space $\mathbb{Y}$ and of its complement $\mathbb{W}$ consisting of "partially complex" 3-planes in $\mathbb{C}^3$.
\end{abstract}

\maketitle
\par\ \par\

\section{Introduction}
\par\ \par\
The subspace $\mathbb{W}$ of "partially complex" 3-planes in the Grassmanian $\textbf{\emph{G}}_{6,3}$ and its complement $\mathbb{Y}$ consisting of totally real 3-planes are important to the study complex tangents and CR-structures. The image of the Gauss map of an embedding (or immersion) of a 3-manifold $M$ into $\mathbb{C}^3$ and its intersection with these complementary subspaces of the Grassmanian determine the complex tangent structure. In this dimension 3, CR-structures are either totally real or totally complex, two extremes which have been subject to much study. In his paper [6], Forsteric resolved the totally real problem in dimension 3; in particular he showed that every closed oriented 3-manifold can be embedded totally real into $\mathbb{C}^3$. In [7], Gromov used the h-principle to demonstrate that the only spheres $S^n$ admitting totally real embeddings into $\mathbb{C}^n$ are $S^1$ and $S^3$.
\par\
The more general situations (proper subsets of complex tangents) have not been as thoroughly investigated. The work of Ahern and Rudin (in [1]) analyzed means of determining complex tangents for specific examples and Webster (in [15]) derived topological invariants for the set of complex tangents of embeddings of real $n$-manifolds into $n$-dimensional complex Euclidean space. In the general situation of embeddings of a real manifold into a complex manifold of arbitrary dimensions, the work of Lia in [11] and later of Domrin in [4] demonstrated formulas relating the characteristic classes of the set of complex tangents with those of the embedded manifold and the ambient complex manifold. Their results hold only under special assumptions on the dimensions of the manifolds and the structure of the complex tangents. Furthermore, Slapar in [14] considered the situation of a generic closed oriented manifold embedded into a complex manifold with codimension 2. In this situation, complex tangents generically are discrete (and finite), and can be classified as either elliptic or hyperbolic by comparing the orientation of the tangent space of the embedded manifold at the complex tangent (which is necessarily complex) with the induced orientation of the tangent space as a complex subspace of the tangent space of the complex ambient manifold. He demonstrated that every pair of complex tangents consisting of one elliptic point and one hyperbolic point can be canceled out; in particular there exists an arbitrarily small isotopy of embeddings so that the set of complex tangents of the isotoped embedding consists exactly of the complex tangents of the original embedding minus the given pair of complex tangents (one elliptic and one hyperbolic). 
\par\
Our situation of interest is that of a closed oriented 3-manifold embedded in $\mathbb{C}^3$. In this setting, complex tangents generically arise as knots (links). In our paper [5] we proved that for every knot type in $S^3$ there exists an embedding $S^3 \hookrightarrow \mathbb{C}^3$ which assumes complex tangents exactly along a knot of the prescribed type (along with one degenerate point). Our work could be readily extended to include all types of links.
\par\
In this paper we consider the non-generic situation where isolated (necessarily degenerate) complex tangents arise. We begin by determining the first and second homotopy groups of the Grassmann submanifold $\mathbb{Y}$. In the process, we derive a topological obstruction to the removal of an isolated (degenerate) complex tangent without affecting the embedding outside of a (small) neighborhood of the point. We then prove that the vanishing of this obstruction is a necessary and sufficient condition for the (local) removal of the complex tangent. We make note of further directions of worthy investigations. 
\par\
We finish the article with some further computations in homotopy and homology for the spaces $\mathbb{W}$ and $\mathbb{Y}$.
\par\ \par\
\section{The Homotopy Groups of $\mathbb{Y}$}
\par\ \par\
In this section, we compute the first two homotopy groups of the space $\mathbb{Y}$ consisting of totally real 3-planes in $\textit{\textbf{G}}_{6,3}$. We use basic techniques and results from algebraic topology, which one can readily find in some literature on the subject; in particular, we refer the reader to Hatcher's book (see [8]).  Consider the submanifold: $\mathbb{W} \subset \textbf{\emph{G}}_{6,3}$ consisting of all "partially complex" 3-planes in $\mathbb{C}^3$ ($=\mathbb{R}^6$). More precisely, $\mathbb{W} = \{P \in \textbf{\emph{G}}_{6,3} | P \cap J(P) \neq \{0\} \}$. The fact that $\mathbb{W}$ is a closed 7-dimensional submanifold of this Grassmanian is well-known in the field, but we prove this explicitly in Section 4.
\par\ \par\
Consider then that the complement $\mathbb{Y} = \textbf{\emph{G}}_{6,3} \setminus \mathbb{W} = \{P \in \textbf{\emph{G}}_{6,3} | P \cap J(P) = \{0\} \}$, i.e. $\mathbb{Y}$ is the space of totally real 3-planes. By our remarks above about $\mathbb{W}$, we see that $\mathbb{Y}$ is an open 9-manifold.
\par\
This dimension ($n=3$) is unique in that both these spaces have manifold structures. In higher dimensions the analogous spaces will have stratified and singular structure, and here we limit our considerations to the dimension $n=3$.
\par\ \par\
We will also consider the corresponded spaces of oriented 3-planes, in particular let $\widetilde{\textbf{G}}_{6,3}$ designate the space of real oriented 3-planes in $\mathbb{C}^3$ and analogously $\widetilde{\mathbb{W}} = \{P \in \widetilde{\textbf{G}}_{6,3} | P \cap J(P) \neq \{0\} \}$ and  $\widetilde{\mathbb{Y}}$ its complement.
\par\ \par\
\begin{theorem}: The manifolds $\mathbb{Y}$ and $\widetilde{\mathbb{Y}}$ are open, connected submanifolds of their respective Grassmannian ambient spaces. Furthermore, they have all the same homotopy groups, in particular we compute the first four homotopy groups to be: $\pi_1(\mathbb{Y}) = \mathbb{Z}$, $\pi_2(\mathbb{Y}) = \pi_3(\mathbb{Y}) = \mathbb{Z}_2, \pi_4(\mathbb{Y}) = 0$  and the same for $\widetilde{\mathbb{Y}}$.
      \end{theorem}
 \par\ \par\
 $\textbf{\emph{Proof:}}$
Note that the matrix group $GL_3 (\mathbb{C})$ acts on $\mathbb{Y}$ by multiplication of a given invertible complex matrix by the standard totally real $\mathbb{R}^3 \subset \mathbb{C}^3$, and the stabilizer group of any plane in $\mathbb{Y}$ is the matrix group: $GL_3 (\mathbb{R})$. As such, we get a fibre bundle: $GL_3 (\mathbb{R}) \hookrightarrow GL_3 (\mathbb{C}) \rightarrow \mathbb{Y}$. As the matrix group $GL_3 (\mathbb{C})$ deformation detracts to the unitary group $U(3)$ and $GL_3 (\mathbb{R})$ deformation retracts to the orthogonal group $O(3)$, we get the induced long exact sequence in homotopy:
\par\ \par\
$... \rightarrow \pi_2(U(3)) \rightarrow \pi_2 (\mathbb{Y}) \rightarrow \pi_1 (O(3)) \rightarrow \pi_1 (U(3)) \rightarrow \pi_1 (\mathbb{Y}) \rightarrow \pi_0 (O(3)) \rightarrow \pi_0 (U(3)) \rightarrow \pi_0 (\mathbb{Y}) \rightarrow 0$
\par\ \par\
However, one finds that the first group: $\pi_2(U(3)) =0$, which one can readily compute or use the well-known fact that $\pi_2$ of any Lie Group is zero. Further, $\pi_1 (O(3)) = \mathbb{Z}_2, \pi_1 (U(3)) = \mathbb{Z}, \pi_0 (O(3)) = \mathbb{Z}_2, \pi_0 (U(3)) = 0$.
\par\
So the beginning of the sequence can be written:
\par\
$0 \rightarrow \pi_2 (\mathbb{Y}) \rightarrow \mathbb{Z}_2 \rightarrow \mathbb{Z} \rightarrow ...$.
\par\ \par\
However, as any homomorphism: $\mathbb{Z}_2 \rightarrow \mathbb{Z}$ must be the zero map, we see that our sequence will decompose into three (short) exact sequences:
\par\ \par\
$0 \rightarrow \pi_2 (\mathbb{Y}) \rightarrow \mathbb{Z}_2 \rightarrow 0$
\par\
$0 \rightarrow \mathbb{Z} \rightarrow \pi_1 (\mathbb{Y}) \rightarrow \mathbb{Z}_2 \rightarrow 0$
\par\
$0 \rightarrow \pi_0(\mathbb{Y}) \rightarrow 0$.
\par\ \par\
The first sequence immediately gives that $\pi_2 (\mathbb{Y}) = \mathbb{Z}_2$ and the third sequence implies that $\pi_0(\mathbb{Y}) = 0$.
\par\ \par\
As for the second sequence, we recall that the long exact sequence in homotopy for a fiber bundle: $F \hookrightarrow E \rightarrow B$ arises from the long exact sequence of the pair $(E, F)$ and the fact that the projection map $p:E \rightarrow B$ induces an isomorphism $p_* : \pi_n (E,F) \rightarrow \pi_n (B)$ for each $n$ (assuming B to be path connected).
\par\ \par\
Hence, in our situation we have that: $\pi_1(\mathbb{Y}) \cong \pi_1(U(3), O(3))$ and the second sequence above comes from the short exact sequence:
\par\ \par\
$0 \rightarrow \pi_1 (U(3)) \rightarrow \pi_1 (U(3), O(3)) \rightarrow \pi_0(O(3)) \rightarrow 0$
\par\
Which simplifies to (as above):
\par\
$0 \rightarrow \mathbb{Z} \rightarrow \pi_1 (U(3),O(3)) \rightarrow \mathbb{Z}_2 \rightarrow 0$
\par\ \par\
Elementary algebra now tells us that there are two possibilities: either
\par\
$\pi_1 (U(3),O(3)) = \mathbb{Z} \oplus \mathbb{Z}_2$ and the first (non-trivial) map is identity onto the first coordinate and the second map is just the map: $(a, b) \rightarrow b$, or it must be the case that $\pi_1 (U(3),O(3)) = \mathbb{Z}$, the first map is multiplication by 2 and the second map is reduction modulo 2.
\par\ \par\
The first map is induced by the quotient map: $j:U(3) \rightarrow U(3)/O(3)$, as $(U(3), O(3))$ is a good pair ($O(3) \subset U(3)$ is a submanifold). Now, it is a well-known fact that the inclusion $U(1) \hookrightarrow U(3)$ induces isomorphism on $\pi_1$, where we may include $U(1)$ as the set of (3x3)-matrices with $e^{i\theta}$ along the diagonal, and zeros everywhere else. As $U(1) \cong S^1$ is a circle, this in fact gives us a generator for $\pi_1 (U(3))$... denote this generating circle by $\gamma$. Now, we note that there are exactly two points on this generating circle which are elements of $O(3)$, namely the identity matrix I and its negative -I. As such, the image of the generating circle $j(\gamma)$ will be a "figure 8" type curve, in particular two circles both passing through the base point in $U(3)/O(3)$. Further, by symmetry it is clear that these circles must represent the same homotopy class in $\pi_1 (U(3),O(3))$, as such we find: $j_*([\gamma]) = [j(\gamma)] = 2 \sigma$, for some class $\sigma \in \pi_1 (U(3),O(3))$. Now, as $j_*$ is injective, $\sigma \neq 0$, and hence is given by multiplication by 2 with some non-trivial class. But then assuming $\pi_1 (U(3),O(3)) = \mathbb{Z} \oplus \mathbb{Z}_2$ will lead to contradiction, as $j_*$ must be identity onto the first coordinate, which cannot now be the case as $\sigma$ itself can never be realized.
\par\
Hence, as there were only two possibilities and the first case cannot hold, it must be that the other possibility holds, as such: $\pi_1 (U(3),O(3)) = \mathbb{Z}$, and $j_*$ is multiplication by 2 and the second map is reduction modulo 2.
\par\ \par\
Therefore, we have proven that  $\pi_1(\mathbb{Y}) \cong \mathbb{Z}$, and we found earlier that $\pi_2(\mathbb{Y}) = \mathbb{Z}_2$ and $\mathbb{Y}$ is path-connected.
\par\
As for $\widetilde{\mathbb{Y}}$, it is clearly a double cover of $\mathbb{Y}$ with an oriented totally real plane mapped to the ambient plane "forgetful" of orientation. Further, as $\widetilde{\mathbb{Y}}$ is the complement of $\widetilde{\mathbb{W}}$ in the space of oriented Grassmannian, both of which are connected and $\widetilde{\mathbb{W}}$ being a closed manifold of codimension 2, it is immediate that $\widetilde{\mathbb{Y}}$ is open and path-connected itself.
\par\ \par\
As such, $\pi_1(\widetilde{\mathbb{Y}}) \mapsto \pi_1(\mathbb{Y})$ is an injective map whose image is of index 2 in  $\pi_1(\mathbb{Y}) = \mathbb{Z}$. Hence, it must be the case that $\pi_1(\widetilde{\mathbb{Y}}) =\mathbb{Z}$. Further, as its a covering it must share all higher homotopy groups with $\mathbb{Y}$ as well.
\par\ \par\
We may further consider the later part of the exact sequence:
\par\ \par\
$\pi_4 (U(3)) \rightarrow \pi_4 (\mathbb{Y}) \rightarrow \pi_3(O(3)) \rightarrow \pi_3(U(3)) \rightarrow \pi_3(\mathbb{Y}) \rightarrow \pi_2(O(3))$
\par\ \par\
First, we notice that $\pi_4(U(3)) = \pi_2(O(3)) = 0$.
\par\ 
Further, we know from Bott in [2] that $\pi_3(O(3)) = \pi_3(U(3)) = \mathbb{Z}$, and that in $\pi_3$ this (inclusion) map is given by multiplication by 2. So we have the sequence:
\par\ \par\
$0 \rightarrow  \pi_4(\mathbb{Y}) \rightarrow \mathbb{Z} \rightarrow \mathbb{Z} \rightarrow  \pi_3(\mathbb{Y}) \rightarrow 0$
\par\ \par\
Since the middle map $\mathbb{Z} \rightarrow \mathbb{Z}$ is multiplication by 2, we immediately get that:
\par\
$\pi_3(\mathbb{Y}) = \mathbb{Z}_2, \pi_4(\mathbb{Y}) = 0$
\par\ \par\
$\textbf{\emph{QED}}$
\par\ \par\
We note to that our techniques above can be continued to find some higher homotopy groups of $\mathbb{Y}$, in particular: $\pi_5(\mathbb{Y}) = \mathbb{Z}, \pi_6(\mathbb{Y}) = \mathbb{Z}_2$. We leave it to the interested reader to verify and work out more of these homotopy groups For the purposes of this paper we will need only the homotopy groups in the above theorem.
\par \par\
\section{Obstruction to the Removal of a Degenerate Complex Tangent}
\par\ \par\
Let $M \hookrightarrow \mathbb{C}^3$ be an embedding, and let $x \in M$ be an isolated degenerate complex tangent. Take a neighborhood of $x$ in $M$, call it $\textbf{B}$, so that $\textbf{B}$ is contained in the interior of an open set whose closure contains no complex tangent other than x. Consider $\textbf{S} = \partial \textbf{B}$ its boundary... as such $\textbf{S}$ contains no complex tangents. Then $\textbf{B} \cong B^3$ and $\textbf{S} \cong S^2$, in particular we get an embedding $i: S^2 \hookrightarrow M$ whose image $i(S^2) = \textbf{S}$. Let $G:M \rightarrow \textbf{\emph{G}}_{6,3}$ denote the Gauss map of the embedded $M \hookrightarrow \mathbb{C}^3$, and consider the composition: $G \circ i:S^2 \rightarrow \textbf{\emph{G}}_{6,3}$. In fact, as all the points of $\textbf{S}$ are totally real, we get a map: $G \circ i:S^2 \rightarrow \mathbb{Y}$. As such, the sphere $\textbf{S}$ gives us an element of $\pi_2(\mathbb{Y}) = \mathbb{Z}_2$. If the isolated complex tangent x where "removable", i.e. if we could slightly perturb the embedding of $M \hookrightarrow \mathbb{C}^3$ near x so that x is not complex tangent nor any point near x is complex tangent, this sphere $G(\textbf{S})$ must be fillable in $\mathbb{Y}$, i.e. the class $[G(\textbf{S})] \in \pi_2(\mathbb{Y})$ must be trivial.
\par\ \par\
Hence, we have found an obstruction to the removal of an isolated complex tangent, namely the homotopy class in $\mathbb{Y}$ of the image of a sphere bounding a sufficiently small neighborhood of x under the Gauss map $G$. Note this will be independent of choice of neighborhood since any two sufficiently small neighborhoods will necessarily have homotopic boundaries (in $\mathbb{Y}$).
\par\ \par\
We note that as $\pi_2(\mathbb{Y}) = \mathbb{Z}_2$, the situation may be that $[G(\textbf{S})] = 1$. In this case, it will not be possible to (locally) remove the complex tangent without affecting the structure of the embedding outside of a neighborhood of the complex tangent. Hence, the vanishing of this obstruction is necessary to "surgically remove" the degenerate point (leaving the embedding the same outside a small neighborhood).
\par\ \par\
We then ask the question: is the vanishing of this obstruction sufficient to ensure the complex tangent $x$ is removable? In other words, if a sphere $\textbf{S} \subset M$ that bounds a (small) neighborhood of $x$ has the property that $[G(\textbf{S})] = 0 \in \pi_2(\mathbb{Y})$ is it necessarily true that we can perturb $M$ slightly near $x$ to make $x$ totally real, without adding new complex tangents?
\par\ \par\
We find that the answer to this question is in the affirmative:
\begin{theorem}:
 Suppose an embedding of a 3-manifold $f: M \hookrightarrow \mathbb{C}^3$ (of class  $\mathcal{C}^k$) with Gauss map $G$ admits an isolated complex tangent $x \in M$ for which there is a neighborhood $\textbf{B}$ of $x$ whose boundary sphere $\textbf{S} \subset M$ satisfies: $[G(\textbf{S})] = 0 \in \pi_2(\mathbb{Y})$. Then there exists a $\mathcal{C}^k$-embedding $\widetilde{f}: M \hookrightarrow \mathbb{C}^3$ that can be taken to be $\mathcal{C}^0$-close to $f$ and regularly homotopic to $f$ so that:
 \par\
 1. $f = \widetilde{f}$ on $M \setminus \textbf{B}$
 \par\
 2. The set of complex tangents of $\widetilde{f}$ equals the set of complex tangents of $f$ minus the point $x$, i.e.:
 $\aleph_{\widetilde{f}} = \aleph_{f} \setminus \{x\}$
      \end{theorem}
\par\ \par\
 $\textbf{\emph{Proof:}}$
 Consider the neighborhood $\textbf{B}$ of $x$ and its boundary sphere $\textbf{S}$ satisfying: $[G(\textbf{S})] = 0 \in \pi_2(\mathbb{Y})$. 
 \par\
 Take an arbitrarily small "thickening" of $\textbf{B}$, call it $\widetilde{\textbf{B}}$. Then $\widetilde{\textbf{B}}$ is a 3-ball in $M$ and $\textbf{A} = \widetilde{\textbf{B}} \setminus \bar{\textbf{B}}$ is an open annulus.
 \par\ \par\
 Now, as  $[G(\textbf{S})] = 0 \in \pi_2(\mathbb{Y})$, we have that $G_* ([S^2]) = 0$ and so there exists a homotopy $H: S^2 \times I \rightarrow \mathbb{Y}$ with:
 \par\
 $H_1 : S^2 \rightarrow \mathbb{Y}$ is $G$ and  $H_0 : S^2 \rightarrow \mathbb{Y}$ is the constant map to a point $*$ in  $\mathbb{Y}$.
 \par\ \par\
 We may now use the homotopy H to define a map: $\widetilde{G}: \bar{\textbf{B}} \rightarrow \mathbb{Y}$ as follows using spherical coordinates $(\theta,\phi,\rho)$:
 \par\
 $\widetilde{G}(\theta,\phi,1) = G(\theta,\phi)$  (on the boundary $S^2$)
 \par\
 $\widetilde{G}(\theta,\phi,\rho) = H_{\rho} (\theta,\phi)$ 
 \par\
 $\widetilde{G}(0,0,0) = *$ 
 \par\ \par\
 By construction and the fact that $H$ is a homotopy, we may assume that this mapping $\widetilde{G}$ is at least continuous. We may further extend this map (using the Gauss map of $f$ on $\textbf{A}$) to obtain a mapping:  
 \par\
 $\widetilde{G}: \widetilde{\textbf{B}} \rightarrow \mathbb{Y}$ that is continuous.
\par\ \par\
Now consider the set of homotopy classes of maps $[B^3,\mathbb{Y};G|_{\textbf{S}}]$ that restrict along the boundary to the map $G$ on $\textbf{S} \cong S^2$. As this set is non-empty (we constructed such a map $\widetilde{G}$ above), the set is in 1-1 correspondence with $\pi_3 (\mathbb{Y})$ (see Hatcher in [8]). 
\par\ \par\  
Our goal is to demonstrate that the Gauss map $G:\textbf{B} \rightarrow \textbf{\textit{G}}_{3,6}$ is homotopic to a map $\hat{G}:\textbf{B} \rightarrow \mathbb{Y}$ that agrees with $G$ on the boundary sphere $\textbf{S}$, and through such maps (agreeing on $\textbf{S}$). Note that $G \in [B^3,\textbf{\textit{G}}_{3,6};G|_{\textbf{S}}] \cong \pi_3(\textbf{\textit{G}}_{3,6})$
\par\ \par\
Let $V_{3,6}$ be the generalized Stiefel manifold of all 3-frames in $\mathbb{R}^6$ and consider the subset of totally real 3-frames ${V_{3,6}}^{tr}$. We get the natural maps sending a frame to the subspace spanned by the frames in $\mathbb{R}^6$:
\par\
$s: V_{3,6} \rightarrow \textbf{\textit{G}}_{3,6}$
\par\
$r: {V_{3,6}}^{tr} \rightarrow \mathbb{Y}$ 
\par\
(see Forstneric in [6] for reference for our work here and below)
\par\ \par\
Notice then the commutative diagram:
\par\ \par\ 
  $\begin{CD}
{V_{3,6}}^{tr} @>>{i_V}>            V_{3,6} \\
 @VVrV @VVsV \\
\mathbb{Y}  @>>{i_G}>          \textbf{\textit{G}}_{3,6}
\end{CD}$
\par\ \par\
where $i_V,i_G$ are inclusion maps.
\par\ \par\
Applying $\pi_3$, we get the resulting commutative diagram of groups:
\par\ \par\
 $\begin{CD}
\pi_3({V_{3,6}}^{tr}) @>>{{i_V}}_*>            \pi_3(V_{3,6}) \\
 @VVr_*V @VVs_*V \\
\pi_3(\mathbb{Y})  @>>{{i_G}_*}>          \pi_3(\textbf{\textit{G}}_{3,6})
\end{CD}$
\par\ \par\
Now, it is known that $\pi_3(V_{3,6}) = \mathbb{Z}_2$, and from the long exact sequence of the fiber bundle over $\textbf{\textit{G}}_{3,6}$, we get that the map $s_*$ is an isomorphism.
\par\ \par\
Also, notice that ${V_{3,6}}^{tr} \cong Gl_3(\mathbb{C})$ retracts to $U(3)$, and that the space of orthonormal 3-frames ${V_{3,6}}^O$ is a retract of $V_{3,6}$. Hence, the map $i_V: {V_{3,6}}^{tr} \rightarrow V_{3,6}$ is homotopy equivalent to the map $\alpha: U(3) \rightarrow {V_{3,6}}^O$, which is the composition of the inclusion $U(3) \hookrightarrow O(6)$ with the quotient map $O(6) \rightarrow {V_{3,6}}^O$. Since the first map is an isomorphism on $\pi_3$ and the second is onto on $\pi_3$ (see [2]), we see that $\alpha_*:\pi_3(U(3)) \rightarrow \pi_3({V_{3,6}}^O)$ is onto, and hence the inclusion map $(i_V)_*$ is an onto map on $\pi_3$.
\par\ \par\
Hence, as the above diagram in $\pi_3$ commutes, $s_* \circ {i_V}_* = {i_G}_* \circ r_*$ is onto, and since $\pi_3(Y) = \mathbb{Z}_2 = \pi_3(\textbf{\textit{G}}_{3,6})$ (as shown above), we have that:
\par\
${i_G}_*: \pi_3(\mathbb{Y})  \rightarrow \pi_3(\textbf{\textit{G}}_{3,6})$ is an isomorphism.
\par\ \par\
Recall that: 
\par\
$\pi_3(\textbf{\textit{G}}_{3,6}) \cong [B^3,\textbf{\textit{G}}_{3,6};G|_{\textbf{S}}]$ and $\pi_3(\mathbb{Y}) \cong [B^3,\mathbb{Y};G|_{\textbf{S}}]$.
\par\ \par\
Now, the Gauss map $G$ of our original embedding is clearly an element of $[B^3,\textbf{\textit{G}}_{3,6};G|_{\textbf{S}}]$, and by the above equivalence and the fact that $i_G$ is the inclusion map, there exists a map $\hat{G} \in [B^3,\mathbb{Y};G|_{\textbf{S}}]$ which is homotopic to $G$ through maps $G_t$ which all agree on the boundary $S^2 \cong \textbf{S}$.
\par\ \par\
Consider then the maps: $\mathcal{G}_t: M \rightarrow \textbf{\textit{G}}_{3,6}$, where ${\mathcal{G}_t}|_{\textbf{B}} = G_t,  {\mathcal{G}_t}|_{M \setminus \textbf{B}} = G$. Then ${\mathcal{G}}_0 = G$ and ${\mathcal{G}_1}|_{\textbf{B}} = \hat{G}$, which is totally real on $\textbf{B}$. All these maps can be taken to be continuous. We may further extend all these maps continuously to the "thickened" ball $\widetilde{B}$ using the original Gauss map $f$.
\par\ \par\
In fact, by the Whitney Approximation Theorem, we may assume that all these maps are of class $\mathcal{C}^k$ (by taking a small homotopy). Furthermore, we may make use the relative version of the Whitney Approximation Theorem  and assume without loss of generality that $\mathcal{G}_t = G$ on the annulus $\textbf{A}$. (See Lee in [10] for an exposition on the Whitney Approximation Theorem).
\par\ \par\
Hence, we have constructed a formal totally real map, which we now denote $\hat{G} = \mathcal{G}_1:\widetilde{\textbf{B}} \rightarrow \mathbb{Y}$ of class  $\mathcal{C}^k$ that agrees with the Gauss map $G$ of the original embedding on a neighborhood $\textbf{A}$ of the boundary, and which is homotopic to $G$ through maps $\widetilde{\textbf{B}} \rightarrow \textbf{\textit{G}}_{3,6}$, In the language of the literature, we have shown that our (original) embedding $f$ is $\textit{formally totally real}$.
\par\ \par\
Let us now consider the $h$-principle for the totally real embeddings. In particular, we have now that our embedding $f:\widetilde{\textbf{B}} \hookrightarrow \mathbb{C}^3$ has Gauss map $G$ which is homotopic (through maps into the Grassmannian) to a map $\hat{G}:\widetilde{\textbf{B}} \rightarrow \mathbb{Y}$, and the restriction of the homotopy on the annulus $\textbf{A} = \widetilde{\textbf{B}} \setminus \textbf{B}$ is equal to the original Gauss map $G$. 
\par\ \par\
As we see from Cieliebak-Eliashberg in [3] and Forstneric in [6], the differential relation of totally real embeddings, or equivalently $\mathbb{Y}$-directed embeddings, satisfies the $h$-principle for extensions. Since the above homotopy is fixed to be $G$ on the annulus $\textbf{A}$ and $G$ is holonomic (totally real) on $\textbf{A}$, we then conclude that there exists a totally real embedding $\widetilde{f}: \widetilde{\textbf{B}} \hookrightarrow \mathbb{C}^3$ that agrees with the original embedding $f$ on the annulus $\textbf{A}$.
\par\ \par\
We give below the formal statement of the $h$-principle for extensions of totally real embeddings:
\par\ \par\ \par\
$\textbf{\textit{\texttt{The h-Principle for Extensions of Totally Real Embeddings}}}$
\par\ \par\
Let $f:M^n \hookrightarrow \mathbb{C}^n$. The $h$-principle for extensions of $\mathcal{C}^k$ totally real embeddings  from a subset $C' \subset M$ to a subset $C \subset M$ which contains $C'$ claims that for every $\mathcal{C}^k$ formally totally real embedding $\varphi_t : Op(C) \rightarrow \textbf{\textit{G}}_{n,2n}$ with $\varphi_0 = G_f$ (the Gauss map of $f$) and satisfying that $\varphi_t|_{Op(C')} = G_f$ is totally real on $Op(C')$ and $Im(\varphi_1) \subset \mathbb{Y}$,  there exists a $\mathcal{C}^k$-homotopy of $f$ to a $\mathcal{C}^k$ (genuine) totally real embedding $\widetilde{f}$ by a homotopy of formally totally real embeddings $\widetilde{f_t}$ ($0 \leq t \leq 1, \widetilde{f_0} =f)$ that is $C^0$-close to $f$ and so that $\widetilde{f_t}|_{Op(C')} = G_f$ is constant in $t$.
 \par\ \par\ \par\
 Now, to apply this $h$-principle (for extensions) to our problem, consider $C = \widetilde{\textbf{B}}$ and $C'= \textbf{A}$. Then the homotopy $\mathcal{G}_t:\widetilde{\textbf{B}} \rightarrow \textbf{\textit{G}}_{3,6}$ is a formally totally real embedding since $\mathcal{G}_0 = G = G_f$ and $Im(\mathcal{G}_1) \subset \mathbb{Y}$, and which is holonomic on the subset $\mathbf{A}$ ($\mathcal{G}_t|_{\textbf{A}}) = G$, the Gauss map of our original embedding $f$.
\par\ \par\
By the above $h$-principle for extensions, there exists a $\mathcal{C}^k$-map $\widetilde{f}: \widetilde{\textbf{B}} \rightarrow \mathbb{C}^3$ which is totally real and satisfies  $\widetilde{f} |_{\textbf{A}} = f$. 
\par\ \par\
Now, we wish to extend $\widetilde{f}$ to the entire of the manifold $M$. Consider that $\widetilde{f}: \widetilde{\textbf{B}} \rightarrow \mathbb{C}^3$ is an embedding. The only concern is for self-intersections for an extension of $\widetilde{f}$ to all of $M$, which cannot be ensured directly from the $h$-principle due to the other complex tangents of the embedding $f$. However, consider that $\widetilde{f} = f$ on the annulus $\textbf{A} = \widetilde{\textbf{B}}\setminus \textbf{B}$. Since $M$ is compact and $f: M \hookrightarrow\mathbb{C}^3$ is an embedding, there will be a minimum distance from the position of any points in $f(M \setminus \widetilde{\textbf{B}})$ to points of the (inner) ball $f(\textbf{B})$. In particular, consider the number $q = inf(dist(f(x),f(y))| x \in M \setminus \widetilde{\textbf{B}}, y \in \textbf{B})$. Note that $q>0$ since $M$ is compact and $f$ is an embedding, and the sets $M \setminus \widetilde{\textbf{B}}$ and $\textbf{B}$ are separated by the (open) annulus $\textbf{A}$. Now, recall that $\widetilde{f}$ is $C^0$-close to $f$, in other words,  $\widetilde{f}$ can be as close as like to $f$ in the (homotoped set) $\textbf{B}$. Suppose then that $dist(f(x), \widetilde{f}(x))|x \in \textbf{B}) < \frac{q}{2}$. Then it is clear that $dist(f(x),\widetilde{f}(y))| x \in M \setminus \widetilde{\textbf{B}}, y \in \textbf{B}) > \frac{q}{2}$. Therefore, there can be no self-intersection of $\widetilde{f}(\textbf{B})$ and $f(M \setminus \widetilde{\textbf{B}})$. As both these sets do not intersect with $f(\textbf{A}) = \widetilde{f}(\textbf{A})$, we may conclude that the extension (for this "close" choice) of $\widetilde{f}$ to all of $M$ by defining $\widetilde{f}(x) = f(x)$ for all $x \in M \setminus {\textbf{B}}$ will have no self intersections. As this extension $\widetilde{f}$ is clearly an immersion (by construction), we may conclude that $\widetilde{f}:M \hookrightarrow \mathbb{C}^3$ is a $\mathcal{C}^k$-embedding which is homotopic to $f$.
\par\ \par\
As $\widetilde{f}$ is totally real on $\widetilde{\textbf{B}}$ and $f$ was complex tangent in $\widetilde{\textbf{B}}$ only at the point $x$, the assertions as claimed in the theorem are satisfied by this choice of $\widetilde{f}$.
\par\ \par\
$\textbf{\textit{QED}}$ 
\par\ \par\
It is of interest to investigate if there are topological conditions on the manifold $M$ that would ensure the vanishing of this invariant. In our paper [5], we found families of embeddings of $S^3$ taking complex tangents along any given class of knot or link which are non-degenerate, union an isolated degenerate point. In some cases the degeneracy was removable by a perturbation, but other cases were not clear. It would be worthwhile to investigate if it were always possible to remove the degenerate complex tangent point while leaving the knot of complex tangents unaffected.
\par\ \par\
We also proved that we can embed $S^3$ so that the complex tangents form exactly a knot of the prescribed type, with one point on the knot being degenerate. Taking a neighborhood containing the knot, we could ask if it is possible to remove the complex tangents while keeping the embedding unchanged outside the neighborhood to obtain a totally real embedding which agrees with the original embedding outside the neighborhood. 
\par\ \par\
In direct analogy with our arguments for the above theorem, we can use the $h$-principle to see that this procedure to (locally) remove the knot of complex tangents would be inherently dependent on the vanishing of the push-forward of the class of the boundary of a neighborhood of the complex tangents into $\mathbb{Y}$. 
\par\ \par\
We believe these are directions of worthy investigation and see the possibility for meaningful and applicable generalizations.
\par\ \par\
\section{Further Homotopy and Homology Groups of $\mathbb{W}$ and $\mathbb{Y}$} 
\par\ \par\
In this section, we will investigate the spaces $\mathbb{W}$ and $\mathbb{Y}$ for themselves and compute some homotopy and homology groups for these interesting spaces.
\par\
\begin{theorem}: The subset $\mathbb{W} \subset \textbf{\emph{G}}_{6,3}$ is a closed orientable 7-submanifold with homotopy groups: $\pi_1 (\mathbb{W}) = \mathbb{Z}_2, \pi_2 (\mathbb{W}) = \pi_3(\mathbb{W}) = \mathbb{Z}$.  The manifold $\widetilde{\mathbb{W}}$ is the universal (double) covering of $\mathbb{W}$, and as such is simply connected and shares all its higher homotopy groups with $\mathbb{W}$.
      \end{theorem}
 \par\
 $\textbf{\emph{Proof:}}$
 We first note that $\mathbb{W}$ is a fiber bundle over $\mathbb{CP}^2$ fibered by real projective space $\mathbb{RP}^3$, i.e. we have a fibre bundle: $\mathbb{RP}^3 \hookrightarrow \mathbb{W} \rightarrow \mathbb{CP}^2$. In particular, given any complex line in $\mathbb{C}^3$, i.e. an element of $\mathbb{CP}^2$, we need to choose a real line orthogonal to it to form an element of $\mathbb{W}$. As a complex line is homeomorphic to $\mathbb{R}^2$ and it is contained in $\mathbb{C}^3 = \mathbb{R}^6$, its complement will be a copy of $\mathbb{R}^4$. Hence, we choose any real line in $\mathbb{R}^4$, i.e. an element of $\mathbb{RP}^3$.
 \par\
 As such, we immediately get that $\mathbb{W}$ is compact and connected (as the base and fibers are). Further, it is orientable as the both the base and the fiber are and the base is simply connected.
 \par\ \par\
 We recall that any fiber bundle (or more generally any fibration) admits a long exact sequence in homotopy... we refer the reader to Hatcher in [8] for its construction. This sequence will take the form:
 \par\ \par\
 $... \rightarrow \pi_4 (\mathbb{CP}^2) \rightarrow \pi_3 (\mathbb{RP}^3) \rightarrow \pi_3 (\mathbb{W}) \rightarrow \pi_3 (\mathbb{CP}^2) \rightarrow \pi_2 (\mathbb{RP}^3) \rightarrow \pi_2 (\mathbb{W}) \rightarrow \pi_2 (\mathbb{CP}^2) \rightarrow \pi_1 (\mathbb{RP}^3) \rightarrow \pi_1 (\mathbb{W}) \rightarrow \pi_1 (\mathbb{CP}^2) \rightarrow \pi_0 (\mathbb{RP}^3) \rightarrow \pi_0 (\mathbb{W}) \rightarrow \pi_0 (\mathbb{CP}^2) \rightarrow 0$
\par\ \par\
Now, first note that both $\mathbb{CP}^2, \mathbb{RP}^3$ are path connected, and as such their $\pi_0 =0$. Further, using the fibration for $\mathbb{CP}^2$ given by: $S^1 \hookrightarrow S^5 \rightarrow \mathbb{CP}^2$, we find that $\mathbb{CP}^2$ is simply connected, $\pi_2(\mathbb{CP}^2) = \mathbb{Z}$, and $\pi_k(\mathbb{CP}^2) = \pi_k(S^5)$, for all $k \geq 3$. In particular, $\pi_3(\mathbb{CP}^2) = \pi_4(\mathbb{CP}^2) = 0$.
\par\
Further, using the fact that $\mathbb{RP}^3$ is double covered by the sphere $S^3$, we find that: $\pi_3(\mathbb{RP}^3) = \mathbb{Z}, \pi_2(\mathbb{RP}^3) = 0, \pi_1(\mathbb{RP}^3) = \mathbb{Z}_2$.
\par\ \par\
As such, our above sequence degenerates to the following (short) exact sequences:
\par\ \par\
$0 \rightarrow \mathbb{Z} \rightarrow \pi_3(\mathbb{W}) \rightarrow 0$
\par\
$0 \rightarrow \pi_2 (\mathbb{W}) \rightarrow \mathbb{Z} \rightarrow \mathbb{Z}_2 \rightarrow \pi_1 (\mathbb{W}) \rightarrow 0$
\par\
$0 \rightarrow \pi_0(\mathbb{W}) \rightarrow 0$.
\par\ \par\
We immediately get that: $\pi_3(\mathbb{W}) = \mathbb{Z}, \pi_0(\mathbb{W}) = 0$ from the first and third sequences.
\par\
As for the second sequence, there are exactly two possibilities for the third map $\mathbb{Z} \rightarrow \mathbb{Z}_2$; it is either onto or the zero map.
\par\
If it were the zero map, then the sequence would split: $0 \rightarrow \pi_2 (\mathbb{W}) \rightarrow \mathbb{Z} \rightarrow 0$ and $0 \rightarrow \mathbb{Z}_2 \rightarrow \pi_1 (\mathbb{W}) \rightarrow 0$, which would imply that $\pi_2 (\mathbb{W}) = \mathbb{Z}, \pi_1 (\mathbb{W}) = \mathbb{Z}_2$.
\par\ \par\
Now, if it were onto, it would necessitate that $\pi_1 (\mathbb{W}) = 0$, and we would get the exact sequence: $0 \rightarrow \pi_2 (\mathbb{W}) \rightarrow \mathbb{Z} \rightarrow \mathbb{Z}_2 \rightarrow 0$, which again by elementary algebra would imply that $\pi_2 (\mathbb{W}) = \mathbb{Z}$ and the second map would be multiplication by two.
\par\ \par\
As these are the only two possibilities, we find that: $\pi_2 (\mathbb{W}) = \mathbb{Z}$, and $\pi_1 (\mathbb{W}) = \mathbb{Z}_2$ or $0$.
\par\ \par\
Hence, we computed $\pi_3, \pi_2, \pi_0$ of $\mathbb{W}$, but we are not sure about the fundamental group. Let's now turn to $\widetilde{\mathbb{W}}$ which will shed some light on the fundamental group of $\mathbb{W}$
\par\
  We may similarly construct a fiber bundle structure for $\widetilde{\mathbb{W}}$, as an oriented 3-plane will contain a complex line in $\mathbb{C}^3$ carrying a natural orientation (from the complex structure), and hence all that is left is to choose an oriented line in the complement of the complex line, which is again $\mathbb{R}^4$. But oriented lines are in one-to-one correspondence with the unit sphere, i.e. $S^3$.
\par\
 Hence, $\widetilde{\mathbb{W}}$ takes the structure of a fiber bundle over $\mathbb{CP}^2$ fibered by spheres $S^3$.
 We find immediately from the exact sequence that:
 \par\ \par\
$... \rightarrow \pi_1 (S^3) \rightarrow \pi_1 (\widetilde{\mathbb{W}}) \rightarrow \pi_1 (\mathbb{CP}^2) \rightarrow ...$.
\par\ \par\
But both $S^3$ and $\mathbb{CP}^2$ are simply connected, which implies $\widetilde{\mathbb{W}}$ is simply connected as well.
\par\
Furthermore, $\widetilde{\mathbb{W}}$ is tautologically a double cover over $\mathbb{W}$, where an oriented 3-plane is sent to itself, "forgetting" orientation. Also, it is a connected manifold as its base and fiber are. Hence, $\widetilde{\mathbb{W}}$ is the universal cover of  $\mathbb{W}$. As such, $\pi_1 (\widetilde{\mathbb{W}})$ must be an index two subgroup of $\pi_1 (\mathbb{W})$. But $\pi_1 (\widetilde{\mathbb{W}})=0$, and as such has order one, which implies that $\pi_1 (\mathbb{W})$ must have order two.
\par\
Hence, $\pi_1 (\mathbb{W}) = \mathbb{Z}_2$.
\par\
$\textbf{\emph{QED}}$
\par\ \par\
One could next ask regarding the homology and cohomology of $\widetilde{\mathbb{W}}$ and $\mathbb{W}$. We were able to prove that:
\par\
 \begin{theorem}: The (integral) cohomology groups of $\widetilde{\mathbb{W}}$ are given by: $H^0 (\widetilde{\mathbb{W}}) = H^2 (\widetilde{\mathbb{W}}) = H^5 (\widetilde{\mathbb{W}}) = H^7 (\widetilde{\mathbb{W}}) = \mathbb{Z}$, and all remaining groups are zero. The analogous result for homology groups holds as well.
 Furthermore, the cohomology ring of $\mathbb{W}$ with $\mathbb{Z}_2$ coefficients is given by: $H^* (\mathbb{W}; \mathbb{Z}_2) = \mathbb{Z} [u, v] / \{u^3 =0, v^4 = uv^2 + u^2\}$
      \end{theorem}
 \par\ \par\
 $\textbf{\emph{Proof:}}$
First let's recall the Gysin sequence for sphere bundles: let $S^{n-1} \hookrightarrow E \rightarrow B$  be a fiber bundle which satisfies a certain orientability condition which always holds when B is simply connected (or we take homology with $\mathbb{Z}_2$ coefficients). Then there is an exact sequence, which we call the Gysin sequence:
\par\ \par\
$... \rightarrow H^{i-n} (B;R) \rightarrow H^{i} (B;R) \rightarrow H^{i} (E;R) \rightarrow H^{i-n+1} (B;R) \rightarrow ...$
\par\ \par\
where the second map $\rightarrow H^{i-n} (B;R) \rightarrow H^{i} (B;R)$ is given by cup product with a certain "Euler class" in $H^{n} (B;R)$, and the following map being $p^*$, where $p:E \rightarrow B$ is the bundle projection map.
\par\
As $H^{i} (B;R) = 0$ for $i<0$ (by definition), the initial part of the Gysin sequence gives us isomorphisms: $p^*:H^{i} (B;R) \cong H^{i} (E;R)$ for $i<n-1$, and the sequence thus begins as:
\par\ \par\
$0 \rightarrow H^{n-1} (B;R) \rightarrow H^{n-1} (E;R) \rightarrow H^{0} (B;R) \rightarrow H^{n} (B;R) \rightarrow H^{n} (E;R) \rightarrow ...$.
\par\ \par\
Here we are interested in the bundle: $S^3 \hookrightarrow \widetilde{\mathbb{W}} \rightarrow \mathbb{CP}^2$ which we now call $\xi$, which satisfies the required orientability condition as $\mathbb{CP}^2$ is simply connected.
\par\
Further, we recall the cohomology of $\mathbb{CP}^2$ (integral coefficients) with : $H^0(\mathbb{CP}^2) = H^2(\mathbb{CP}^2)=H^4(\mathbb{CP}^2)=\mathbb{Z}$, and all other cohomology groups are zero. Further, the cohomology ring is generated by the generator $u \in H^2(\mathbb{CP}^2)$, i.e. $u^2 \in H^4(\mathbb{CP}^2)$ generates its group as well.
\par\ \par\
Hence, applying the Gysin sequence in low dimensions, we immediately find:
\par\
$H^0(\widetilde{\mathbb{W}}) = H^2(\widetilde{\mathbb{W}})= \mathbb{Z}$, and $H^1(\widetilde{\mathbb{W}}) =0$. For the remaining cohomology groups, we must investigate the sequence:
\par\ \par\
$0 \rightarrow H^3(\mathbb{CP}^2) \rightarrow H^3(\widetilde{\mathbb{W}}) \rightarrow H^0(\mathbb{CP}^2) \rightarrow H^4(\mathbb{CP}^2) \rightarrow H^4(\widetilde{\mathbb{W}}) \rightarrow H^1(\mathbb{CP}^2) \rightarrow H^5(\mathbb{CP}^2) \rightarrow H^5(\widetilde{\mathbb{W}}) \rightarrow H^2(\mathbb{CP}^2) \rightarrow H^6(\mathbb{CP}^2) \rightarrow H^6(\widetilde{\mathbb{W}}) \rightarrow H^3(\mathbb{CP}^2) \rightarrow H^7(\mathbb{CP}^2) \rightarrow H^7(\widetilde{\mathbb{W}}) \rightarrow H^4(\mathbb{CP}^2)
\rightarrow H^8(\mathbb{CP}^2) \rightarrow ...$
\par\ \par\
From our above knowledge of the cohomology groups of $\mathbb{CP}^2$, we immediately deduce the following:
\par\
$H^5(\widetilde{\mathbb{W}}) = \mathbb{Z}, H^6(\widetilde{\mathbb{W}}) = 0, H^7(\widetilde{\mathbb{W}}) = \mathbb{Z}$. Further, as $\mathbb{W}$ is a 7-manifold, all the higher groups must be zero.
\par\ \par\
 The beginning of the sequence degenerates to: $0 \rightarrow H^3(\widetilde{\mathbb{W}}) \rightarrow H^0(\mathbb{CP}^2) \rightarrow H^4(\mathbb{CP}^2) \rightarrow H^4(\widetilde{\mathbb{W}}) \rightarrow 0$, and both $H^0(\mathbb{CP}^2) = H^4(\mathbb{CP}^2) = \mathbb{Z}$, and the map connecting them is cup product with the Euler class in $H^4(\mathbb{CP}^2)$. So we get:
 \par\
 $0 \rightarrow H^3(\widetilde{\mathbb{W}}) \rightarrow \mathbb{Z} \rightarrow \mathbb{Z} \rightarrow H^4(\widetilde{\mathbb{W}}) \rightarrow 0$.
 \par\ \par\
To find this Euler class, we first make note that our bundle $\xi$ is the unit sphere bundle of a complex 2-bundle; in particular, let: $\eta$ be the bundle: $\mathbb{C}^2 \hookrightarrow E \rightarrow \mathbb{CP}^2$, where the fiber over a complex line is its orthogonal complement in $\mathbb{C}^3$. By construction, it is immediate that $\xi$ is the unit sphere bundle of $\eta$.
\par\ \par\
It is a result from the theory of characteristic classes that the Euler class of the unit sphere bundle of a complex k-bundle is equal the $k^{th}$ Chern class of the complex bundle; see Milnor and Stasheff in [12].
\par\
Next, consider the tautological line bundle over $\mathbb{CP}^2$, call it $\tau$, given by: $\mathbb{C} \hookrightarrow T \rightarrow \mathbb{CP}^2$ where the fiber over a complex line is the collection of points in the complex line. We note that the direct sum of complex bundles: $\tau \oplus \eta = \epsilon$, where $\epsilon$ is the trivial complex 3-bundle over $\mathbb{CP}^2$. Hence, applying the sum formula for the total Chern class of vector bundles, we find:
\par\ \par\
$c(\tau) \smile c(\eta) = c(\epsilon) =1$, as the Chern classes of the trivial bundle are trivial.
\par\
We know that the first Chern class of the tautological bundle is the generator $u \in H^2(\mathbb{CP}^2)$, and as it is a complex line bundle, it has no other Chern classes. Hence, our formula may be written:
\par\ \par\
$(1+u) \smile (1+c_1+c_2) =1$, where $c_1, c_2$ are the first and second Chern classes of $\eta$. Multiplying out and collecting terms of the same degree:
\par\ \par\
$1 + (u +c_1) + (u c_1 +c_2) = 1$, and as such: $c_1 = -u$. Furthermore:
\par\
$0 = uc_1 +c_2 = -u^2 +c_2$, which implies that $c_2 = u^2$, which is the (positive) generator of $H^2(\mathbb{CP}^2)$.
\par\ \par\
Hence, $e(\xi) = c_2(\eta) = u^2$, and as such cup product with the Euler class induces isomorphism on cohomology.
\par\
Recall our sequence:
\par\ \par\
$0 \rightarrow H^3(\widetilde{\mathbb{W}}) \rightarrow \mathbb{Z} \rightarrow \mathbb{Z} \rightarrow H^4(\widetilde{\mathbb{W}}) \rightarrow 0$.
\par\ \par\
As the middle map $ \mathbb{Z} \rightarrow \mathbb{Z}$ is given by cup product with the Euler class and thus an isomorphism, it must be the case that the outer groups are zero (by exactness), in particular: $H^3(\widetilde{\mathbb{W}}) = H^4(\widetilde{\mathbb{W}}) = 0$.
\par\ \par\
 Hence, we have computed all the cohomology groups of $\widetilde{\mathbb{W}}$ with integral coefficients, as $\widetilde{\mathbb{W}}$ is a 7-manifold and all higher groups are zero. Further, as $\widetilde{\mathbb{W}}$ is closed and orientable,
we may apply Poincare duality to achieve all the homology groups of $\widetilde{\mathbb{W}}$, which turn out to be exactly the same in the same dimensions.
\par\ \par\
Now, we recall the following result, for which we refer the reader to Hatcher in [9]. Given a real k-dimensional vector bundle $\eta: E(\eta) \rightarrow X$, we may find the cohomology ring with $\mathbb{Z}_2$-coefficients of the unit sphere bundle with antipodal points in each fiber attached, which we may write as: $S(\eta)/ \mathbb{Z}_2$ in the following way:
\par\ \par\
$H^*(S(\eta)/ \mathbb{Z}_2; \mathbb{Z}_2) = H^*(X; \mathbb{Z}_2)[v] / \{v^k = \Sigma w_i(\eta) v^{k-i}\}$, where the sum in the given relation is for $ 0 < i \leq k$, and $w_i (\eta)$ is the $i^{th}$ Steifel-Whitney class of $\eta$.
\par\ \par\
Applying this result to our situation, consider $\eta$ to be the complex 2-bundle: $\mathbb{C}^2 \hookrightarrow E \rightarrow \mathbb{CP}^2$ as we used above. Then the unit sphere bundle of $\eta$ will be the familiar sphere bundle: $\xi: S^3 \hookrightarrow \widetilde{\mathbb{W}} \rightarrow \mathbb{CP}^2$. Further, by construction $\xi/ \mathbb{Z}_2$ identifying antipodal points in each fiber will give us the bundle $\mathbb{RP}^3 \hookrightarrow \mathbb{W} \rightarrow \mathbb{CP}^2$, whose total space is $\mathbb{W}$.
\par\
Now, considering $\eta$ as a real 4-bundle, we compute the Steifel-Whitney classes to be the reduction modulo 2 of the Chern classes in the even dimensions 2 and 4, and zero otherwise. As such, $w_2(\eta) = u \in H^2(\mathbb{CP}^2; \mathbb{Z}_2), w_4 (\eta) = u^2 \in H^4(\mathbb{CP}^2; \mathbb{Z}_2)$, both generators for their respective groups (both isomorphic to $\mathbb{Z}_2$), and $w_1(\eta)=w_3(\eta) = 0$.
\par\ \par\
In fact, $H^*(\mathbb{CP}^2; \mathbb{Z}_2) = \mathbb{Z}_2 [u]/ \{u^3 = 0\}$ and the relevant relation becomes: $v^4 = u v^2 +u^2$. So we may write:
\par\ \par\
$H^*(\mathbb{W}; \mathbb{Z}_2) = H^*(S(\eta)/ \mathbb{Z}_2; \mathbb{Z}_2) = \mathbb{Z}_2 [u, v]/ \{u^3 = 0, v^4 = u^2 + u v^2\}$, which is what we wished to prove. An easy computation now gives the cohomology groups of $\mathbb{W}$ in  $\mathbb{Z}_2$-coefficients, in particular $\mathbb{Z}_2 \oplus \mathbb{Z}_2$ in dimensions 2, 3, 4, 5, and the group $\mathbb{Z}_2$ in dimensions 0, 1, 6, and 7.
\par\ \par\
$\textbf{\emph{QED}}$
\par\ \par\
We noted above the cohomology groups $H^*(\mathbb{W}; \mathbb{Z}_2)$. Now, consider the projection map $p: \widetilde{\mathbb{W}} \rightarrow \mathbb{W}$, and its induced map on cohomology: $p^*:H^*(\mathbb{W}) \rightarrow H^*(\widetilde{\mathbb{W}})$. We know further that there exists a transfer map $\rho: H^*(\widetilde{\mathbb{W}}) \rightarrow H^*(\mathbb{W})  $ given by the $\mathbb{Z}_2$ action on $\widetilde{\mathbb{W}}$ so that: $\rho \circ p^* : H^*(\mathbb{W})  \rightarrow H^*(\mathbb{W})$ is merely multiplication by 2. Furthermore, with our understanding of $H^*(\widetilde{\mathbb{W}})$ which consists of only $\mathbb{Z}$ summands, this gives us two important facts:
\par\
1.) $H^*(\mathbb{W})$ will consist of only copies (and direct sums) of the groups $\mathbb{Z}$ and $\mathbb{Z}_2$.
\par\
2.) $H^*(\mathbb{W}; \mathbb{Z}[\frac{1}{2}]) \cong H^*(\widetilde{\mathbb{W}}; \mathbb{Z}[\frac{1}{2}])$
\par\
(see [13] for reference)
\par\ \par\
Now, with our knowledge of the groups $H^*(\mathbb{W}; \mathbb{Z}_2)$ and $H^*(\widetilde{\mathbb{W}})$ we may determine the (integral) cohomology groups of $\mathbb{W}$. However, as we may use Poincare Duality, let's first determine the homology groups. Note that: $H^*(\mathbb{W}; \mathbb{Z}_2) \cong H_*(\mathbb{W}; \mathbb{Z}_2)$, which we computed above as being a copy of $\mathbb{Z}_2 \oplus \mathbb{Z}_2$ in dimensions 2, 3, 4, 5, and  $\mathbb{Z}_2$ in dimensions 0, 1, 6, and 7. Also, as $H^*(\mathbb{W}; \mathbb{Z}[\frac{1}{2}]) \cong H^*(\widetilde{\mathbb{W}}; \mathbb{Z}[\frac{1}{2}])$ we know that the only homology groups that admit $\mathbb{Z}$ summands will be 0, 2, 5, and 7... furthermore there must be precisely one such summand in those dimensions. Now, let's recall the Universal Coefficient Theorem for Homology, which states:
\par\ \par\
$0 \rightarrow H_k(X) \otimes G \rightarrow H_k(X; G) \rightarrow Tor(H_{k-1}(X), G) \rightarrow 0$
\par\
are split short exact sequences.
\par\ \par\
First, we get $H_0(\mathbb{W}) = \mathbb{Z}$ as it is connected. Next, since $:H^*(\mathbb{W}; \mathbb{Z}[\frac{1}{2}]) \cong H^*(\widetilde{\mathbb{W}}; \mathbb{Z}[\frac{1}{2}])$, the only dimensions where $\mathbb{W}$ could admit $\mathbb{Z}$ terms in cohomology will be 0, 2, 5, 7 (from cohomology of $\widetilde{\mathbb{W}}$). All other dimensions must consist solely of $\mathbb{Z}_2$ summands (or 0). We get the same result in homology by duality. Now, $H_1(\mathbb{W}) = \mathbb{Z}_2$ since $H_1(\mathbb{W}; \mathbb{Z}_2) = \mathbb{Z}_2$. Further, as $H_2(\mathbb{W}; \mathbb{Z}_2) = \mathbb{Z}_2 \oplus \mathbb{Z}_2$ and we must have a $\mathbb{Z}$ summand, using the UCT we see that $H_2(\mathbb{W}) = \mathbb{Z}$ (we get a $\mathbb{Z}_2$ summand from the previous dimension as well). Now, as $H_2(\mathbb{W})$ is free and $H_3(\mathbb{W})$ has no $\mathbb{Z}$ summand, we see that $H_3(\mathbb{W})= \mathbb{Z}_2 \oplus \mathbb{Z}_2$ from homology with $\mathbb{Z}_2$ coefficients. As this term will arise naturally in the computation of $H_4(\mathbb{W}; \mathbb{Z}_2)$, we see that $H_4(\mathbb{W})$ can have no terms, i.e. must be trivial. Now, $H_5(\mathbb{W})$ admits (exactly) one $\mathbb{Z}$ summand and from the fact that $H_4(\mathbb{W})$ is free (in fact, trivial), the UCT and the above computations with $\mathbb{Z}_2$ coefficients implies that $H_5(\mathbb{W})= \mathbb{Z} \oplus \mathbb{Z}_2$. Also, as $H_6(\mathbb{W}; \mathbb{Z}_2) = \mathbb{Z}_2$ and a $\mathbb{Z}_2$ arises from $H_5(\mathbb{W})$ by the UCT, we see that $H_6(\mathbb{W}) = 0$. Finally, $H_7(\mathbb{W}) = \mathbb{Z}$ as there was no torsion in $H_6(\mathbb{W})$ and the fact that $H_7(\mathbb{W}) $ must admit a $\mathbb{Z}$ summand. Applying Poincare Duality, we summarize our results in cohomology in the following theorem:
\begin{theorem}: $H^0(\mathbb{W}) = H^5(\mathbb{W}) = H^7(\mathbb{W}) = \mathbb{Z}$, $H^1(\mathbb{W}) = H^3(\mathbb{W}) = 0$, $H^2(\mathbb{W}) = \mathbb{Z} \oplus \mathbb{Z}_2$, $H^4(\mathbb{W}) = \mathbb{Z}_2 \oplus \mathbb{Z}_2$, and $H^6(\mathbb{W}) = \mathbb{Z}_2$
      \end{theorem}
 \par\ \par\
Let us now turn to the manifold $\mathbb{Y}$ and its double cover  
$\widetilde{\mathbb{Y}}$.
 First, we compute some homology classes for $\widetilde{\mathbb{Y}}$.
 \begin{theorem}:
 The first two homology groups of $\widetilde{\mathbb{Y}}$ are given by: $H_1(\widetilde{\mathbb{Y}}) = \mathbb{Z}, H_2 (\widetilde{\mathbb{Y}}) = \mathbb{Z}_2$. The first two cohomology groups are given by: $H^1(\widetilde{\mathbb{Y}}) = \mathbb{Z}, H^2(\widetilde{\mathbb{Y}}) = 0$.
      \end{theorem}
 $\textbf{\emph{Proof:}}$ To simplify notation, we write $X= \widetilde{\textbf{G}}_{6,3}$. Take a (small) tubular neighborhood of $\widetilde{\mathbb{W}}$, which we will denote by $V$. Consider the pair: $(X, V)$. Then $H^*(X, V) = H^*(X \setminus \widetilde{\mathbb{W}}, V \setminus \widetilde{\mathbb{W}})$ by excision. This is in turn equals to $H^*(X \setminus \widetilde{\mathbb{W}}, \partial)$ by deformation retraction, where $\partial$ denotes the boundary of the space $X \setminus \widetilde{\mathbb{W}}$.  Now, applying Lefshetz Duality, we find: $H^k(X, \widetilde{\mathbb{W}}) = H^k(X, V) = H^k(X \setminus \widetilde{\mathbb{W}}, \partial) = H_{9-k} (X \setminus \widetilde{\mathbb{W}}) = H_{9-k} (\widetilde{\mathbb{Y}})$.
 \par\ \par\
 Let's look at (a portion of) the long exact sequence of the pair $(X, \widetilde{\mathbb{W}})$ in cohomology:
 \par\ \par\
 $H^6 (\widetilde{\mathbb{W}}) \rightarrow H^7 (X, \widetilde{\mathbb{W}}) \rightarrow H^7(X) \rightarrow H^7 (\widetilde{\mathbb{W}}) \rightarrow H^8 (X, \widetilde{\mathbb{W}}) \rightarrow H^8(X)$
\par\ \par\
First, to find some groups of $X$, we note that the (unoriented) Grassmannian manifold $\textbf{\emph{G}}_{n,k}$ admits a bundle projection $p:V_{n, k} \rightarrow \textbf{\emph{G}}_{n,k}$ fibered by the orthogonal group $O(k)$. Further, as the space $V_{n, k}$ is (n-k-1)-connected, we obtain the following isomorphisms from the long exact sequence of the fiber bundle, for all $i < n-k$:
\par\
 $\pi_i(\textbf{\emph{G}}_{n,k}) = \pi_{i-1} (O(k))$... see [8] for reference.
 \par\ \par\
 Hence, in our situation, we find that: $\pi_1 (\textbf{\emph{G}}_{6,3}) = \pi_2 (\textbf{\emph{G}}_{6,3}) = \mathbb{Z}_2$, from the groups of $O(3)$. \par\
 As $X$ is the universal cover for $\textbf{\emph{G}}_{6,3}$, we find immediately that: $\pi_1 (X) = 0, \pi_2 (X) = \mathbb{Z}_2$, which by Hurewicz Theorem implies that: $H_1(X) = 0, H_2(X) = \mathbb{Z}_2$. Further, by Poincare Duality we find that: $H^7(X) = \mathbb{Z}_2, H^8(X) = 0$.
 \par\ \par\
Together with our knowledge of cohomology of $\widetilde{\mathbb{W}}$, the above sequence may be written:
\par\ \par\
 $0 \rightarrow H^7 (X, \widetilde{\mathbb{W}}) \rightarrow \mathbb{Z}_2 \rightarrow \mathbb{Z} \rightarrow H^8 (X, \widetilde{\mathbb{W}}) \rightarrow 0 $
 \par\
 Further, as any map $\mathbb{Z}_2 \rightarrow \mathbb{Z}$ must be trivial, the sequence breaks up and we immediately get that:
 \par\
 $H_2(\widetilde{\mathbb{Y}}) = \mathbb{Z}_2, H_1 (\widetilde{\mathbb{Y}}) = \mathbb{Z}$, as $H^k(X, \widetilde{\mathbb{W}}) = H_{9-k} (\widetilde{\mathbb{Y}})$
 \par\ \par\
Furthermore, $H^1(\widetilde{\mathbb{Y}}) = Hom(H_1 (\widetilde{\mathbb{Y}}, \mathbb{Z}) = \mathbb{Z}$ and since $ H_1 (\widetilde{\mathbb{Y}}) $ is free, we obtain that: $H^2(\widetilde{\mathbb{Y}}) = Hom(\mathbb{Z}_2, \mathbb{Z}) = 0$, as again any homomorphism $\mathbb{Z}_2 \rightarrow \mathbb{Z}$ must be trivial.
\par\ \par\
$\textbf{\emph{QED}}$
\par\ \par\
We believe these computations in homology and homotopy for the spaces $\mathbb{W}$ and $\mathbb{Y}$ could prove to useful in future research in the field of real submanifolds of complex space.
\par\ 
\section*{Acknowledgement}
The author was supported in part by IRG grant no. 4010501011418 from Alfaisal University.
\par\

\end{document}